# Forbidden Subgraph Characterization of Quasi-line Graphs

Medha Dhurandhar
mdhurandhar@gmail.com

**Abstract:**
Here in particular, we give a characterization of Quasi-line Graphs in terms of forbidden induced subgraphs. In general, we prove a necessary and sufficient condition for a graph to be a union of two cliques.

**1. Introduction:** A graph is a quasi-line graph if for every vertex v, the set of neighbours of v is expressible as the union of two cliques. Such graphs are more general than line graphs, but less general than claw-free graphs. In [2] Chudnovsky and Seymour gave a constructive characterization of quasi-line graphs. An alternative characterization of quasi-line graphs is given in [3] stating that a graph has a fuzzy reconstruction iff it is a quasi-line graph and also in [4] using the concept of sums of Hoffman graphs. Here we characterize quasi-line graphs in terms of the forbidden induced subgraphs like line graphs.

We consider in this paper only finite, simple, connected, undirected graphs. The vertex set of G is denoted by V(G), the edge set by E(G), the maximum degree of vertices in G by $\Delta(G)$, the maximum clique size by $\omega(G)$ and the chromatic number by $\chi(G)$. N(u) denotes the neighbourhood of u and $\overline{N(u)} = N(u) + u$.

For further notation please refer to Harary [3].

**2. Main Result:**

Before proving the main result we prove some lemmas, which will be used later.

**Lemma 1:** If G is $\{3K_1, C_5\}$-free, then either
1) $G \sim K_{|V(G)|}$ or
2) If $v, w \in V(G)$ are s.t. $vw \notin E(G)$, then $V(G) = \{v, w\} \cup B \cup C \cup A_1 \cup A_2 \cup A_3$ where $<B>$, $<C>$, $<A_i>$ (i = 1, 2) are complete.

Proof: Let $G \neq K_{|V(G)|}$. Define $B = \{b \in V(G) / bv \in E(G)$ and $bw \notin E(G)\}$, $C = \{c \in V(G) / cw \in E(G)$ and $cv \notin E(G)\}$, $A_1 = \{a \in V(G) / av, aw \in E(G)$ and $ac \notin E(G)$ for some $c \in C\}$, $A_2 = \{a \in V(G) / av, aw \in E(G)$ and $ab \notin E(G)$ for some $b \in B\}$ and $A_3 = V(G) - (\{v, w\} \cup B \cup C \cup A_1 \cup A_2)$. As G is $3K_1$-free, $<B>$, $<C>$ are complete. If $<A_1>$ is not complete, then let $a, a' \in A_1$ be s.t. $aa' \notin E(G)$. By definition $\exists c \in C$ s.t. $ac \notin E(G)$. As G is $3K_1$-free, $a'c \in E(G)$ and $\exists c' \in C$ s.t. $a'c' \notin E(G)$. But then $<v, a, c', c, a'> = C_5$, a contradiction. Hence $<A_1>$ is complete. Similarly $<A_2>$ is complete.

This proves the Lemma.

**Theorm:** A graph G is a union of two cliques iff G does not have $C^C_{2n+1}$ as an induced subgraph.

Proof: (=>) Obvious

(<=) Let if possible G be a smallest graph without $C^C_{2n+1}$ as an induced subgraph, which is not a union of two cliques. Let $u \in V(G)$. Then clearly deg $u \leq |V(G)|-3$. Let $G-u = \bigcup_{i=1}^{2} Q_i$ where $Q_i$ is complete for i = 1, 2 and $Q_1$ is a maximal such clique i.e. $\forall v \in Q_2 \exists w \in Q_1$ s.t. $vw \notin E(G)$. Let $v_i \in Q_i$ s.t. $uv_i \notin E(G)\}$, i = 1, 2. Define $H^1_i = \{v \in Q_i / vu \notin E(G)\}$. Then $H^1_i \neq \phi$, i = 1, 2. Also as G is $\{C_3^C = 3K_1\}$-free, $vw \in E(G) \ \forall \ v \in Q_1, w \in Q_2$. Iteratively define $H_i^j = \{v \in Q_i / vx \notin E(G)$ for some $x \in H_k^{j-1}$, where $k \in \{1, 2\}$ and $k \neq i\}$ for $i = 1, 2$.

**Claim:** $\forall \ v \in H_1^j$ and $w \in H_2^j$, $vw \in E(G)$ for $j = 1, ..., n$.

Let if possible $\exists \ v_1^j \in H_1^j$ and $w_1^j \in H_2^j$, $v_1^j w_1^j \notin E(G)$. Let $w_2^{j-1} \in H_2^{j-1}$ and $v_1^{j-1} \in H_1^{j-1}$ be s.t. $v_2^{j-1} w_2^{j-1} \notin E(G)$ and so on. Thus we get $<\bigcup_1^t v, \bigcup_1^t w, u> = C_{2j+1}^C$ where $t = 1,...,j$, a contradiction. Thus the **Claim** holds.

Let $R_1 = \{u \cup H_j^{2i}\}$, $1 \leq i \leq \frac{|V(G)|}{2}$, $j = 1, 2$ and $R_2 = \cup H_j^{2i+1}$, $1 \leq i \leq \frac{|V(G)|-1}{2}$, $j = 1, 2$. Then $G = <R_1> \cup <R_2>$ is the required decomposition.

**Corollary 1: G is a Quasi-line graph iff G does not have $K_1 + C_{2n+1}^C$ as an induced subgraph.**

**Corollary 2: If G does not have any of $K_{1,3}$, $W_6$, $K_1 + \{2K_1 + (K_2 \cup K_1)\}$ as an induced subgraph, then G is a Quasi-line graph.**

Proof: Obvious as $\{2K_1 + (K_2 \cup K_1)\}$ is an induced subgraph of $C_{2n+1}^C \ \forall \ n > 2$.

**Remarks**

**1.** Condition $K_1+$ is necessary in Corollary 1, as $C_{2n+1}^C$ itself is Quasi-line.

**2.** As $K_1 + \{2K_1 + (K_2 \cup K_1)\}$ itself is quasi-line, clearly absence of $K_1 + \{2K_1 + (K_2 \cup K_1)\}$ is not necessary for G to be quasi-line.

**References:**

[1] Harary, Frank, *Graph Theory* (1969), Addison–Wesley, Reading, MA.

[2] "Claw-free Graphs. VII. Quasi-line graphs", Maria Chudnovsky, Paul Seymour, Journal of Combinatorial Theory Series B, Volume 102 Issue 6, Nov. 2012, 1267-1294.

[3] "A characterization of quasi-line graphs", Akihiro Munemasa, Tetsuji Taniguchi, July 10, 2010

[4] "Competition Numbers, Quasi-line Graphs, and Holes", Brendan D. McKay, Pascal Schweitzer, and Patrick Schweitzer, *SIAM J. Discrete Math.*, 28(1), 77–91.